\documentclass[runningheads]{lncse}

\usepackage{epsfig}
\usepackage{graphics} 

\usepackage{amsmath,amsfonts,amssymb}
\usepackage{graphicx,verbatim,mathtools}
\usepackage{listings,color}
\usepackage{tikz,setspace}
\usetikzlibrary{patterns}
\usetikzlibrary{positioning}
\usepackage{pgfplots,bm}
\usepackage[only,llbracket,rrbracket]{stmaryrd}

\newcommand{\bo}[1]{\mathbf{#1}} 
\newcommand{\llb}{\llbracket}
\newcommand{\rrb}{\rrbracket}
\newcommand{\lla}{\left\lbrace}
\newcommand{\rra}{\right\rbrace}
\newcommand{\dd}{\mathrm{d}}
\renewcommand{\dim}{d}
\newcommand{\tends}{\rightarrow}

\newcommand{\norm}[1]{\lVert#1\rVert} 
\newcommand{\p}{\partial}
\newcommand{\half}{\textstyle{\frac{1}{2}}}

\DeclareMathOperator{\diam}{diam}

\DeclareMathOperator{\Div}{div}

\newcommand{\eval}[2]{\left. #1\right|_{#2}}

\newcommand{\R}{\mathbb{R}}

\newcommand{\Om}{\Omega}
\newcommand{\DO}{\partial\Om}

\newcommand{\elm}{K}
\newcommand{\face}{F}
\newcommand{\calFhi}{\calF_h^i} 
\newcommand{\calFhb}{\calF_h^b}
\newcommand{\calFh}{\calF_h}
\newcommand{\Vh}{V_{h,k}} 
\newcommand{\Sh}{\mathbf{\Sigma}_{h,k+1}} 

\newcommand{\calF}{\mathcal{F}}
\newcommand{\calT}{\mathcal{T}}

\newcommand{\calP}{\mathcal{P}}
\newcommand{\bcalP}{\bm{\mathcal{P}}}

\newcommand{\bx}{\bm x}
\newcommand{\btau}{\bm \tau}
\newcommand{\bmu}{\bm \mu}
\newcommand{\bsigma}{\bm \sigma}
\newcommand{\bn}{\bm n}

\newcommand{\LRTN}{\mathbf{RTN}_{k+1}(\elm)} 

\begin{document}

\title{Stable discontinuous Galerkin FEM without penalty parameters}
\titlerunning{Stable DGFEM without penalty parameters}
\author{Lorenz John\inst{1} \and Michael Neilan\inst{2} \and Iain Smears\inst{3}}
\authorrunning{John et al.}   
\institute{Technische Universit\"at M\"unchen, M\"unchen,  80333, Germany, {\tt john@ma.tum.de} \and University of Pittsburgh, Pittsburgh, PA 15260, United States, {\tt neilan@pitt.edu} \and INRIA Paris-Rocquencourt, Le Chesnay, 78153, France, {\tt iain.smears@inria.fr}}

\maketitle

\begin{abstract}
We propose a modified local discontinuous Galerkin (LDG)
method for second--order elliptic problems that does not require extrinsic penalization to ensure stability. Stability is instead achieved by showing a discrete Poincar\'e--Friedrichs inequality for the discrete gradient that employs a lifting of the jumps with one polynomial degree higher than the scalar approximation space. Our analysis covers rather general simplicial meshes with the possibility of hanging nodes.
\end{abstract}

\section{Introduction}\label{sec:intro}

It is well--known that the local discontinuous Galerkin (LDG) method for second--order elliptic problems can be formulated, in part,  by replacing the differential operators in the variational formulation by their discrete counterparts \cite{Castillo2000,Cockburn1998,DiPietro2012}.
For example, on the space of discontinuous piecewise polynomials of degree at most $k$, the discrete gradient operator is composed of the element-wise gradient corrected by a lifting of the jumps into the space of piecewise polynomial vector fields.
The original formulation of the LDG method \cite{Castillo2000} employs liftings of same polynomial degree $k$ as the scalar finite element space, while liftings of order $k-1$ have also been considered, see the textbook \cite{DiPietro2012} and the references therein. Part of the motivation for these choices of the order of the lifting is the correspondence to the order of the element-wise gradient and reasons of ease of implementation.
However, unlike the continuous gradient acting on the space $H^1_0$, the discrete gradient operators with liftings of order $k-1$ or $k$ fail to satisfy a discrete Poincar\'e--Friedrichs inequality.
Therefore, the LDG method requires additional penalization with user--defined penalty parameters to ensure stability.

In this note, we construct a modified LDG method with guaranteed stability without the need for extrinsic penalization. This result is obtained by simply increasing the polynomial degree of the lifting operator to order $k+1$ and exploiting properties of the piecewise Raviart--Thomas--N\'ed\'elec finite element space. Our analysis covers the case of meshes with hanging nodes under a mild condition of \emph{face regularity} which we introduce in this work. We recall that the order of the lifting in the LDG method does not alter the dimension or stencil of the resulting stiffness matrix. As a result, the proposed method has a negligible increase of computational cost and inherits the advantages of the standard LDG method in terms of locality and conservativity.

The rest of the paper is organized
as follows.  In Sec.~\ref{sec:notation} we
give the notation used throughout the manuscript
and state some preliminary results.  We define
the lifted gradient operator with increased
polynomial degree in Sec.~\ref{sec:stability}
and show that the $L^2$ norm of this operator
is equivalent to a discrete $H^1$ norm on piecewise polynomial
spaces. We establish by means of a counterexample
that the increased polynomial degree is necessary
to obtain this stability estimate in Sec.~\ref{sec:counter}.
In Sec.~\ref{sec:LDG} we propose and study the modified LDG method in the context of the Poisson equation.

\section{Notation}\label{sec:notation}
\begin{figure}
\begin{center}
\begin{tikzpicture}[scale=1.5,
st1/.style={circle,draw=black,fill=black,thick, minimum size=1mm, inner sep=0pt}
]
\draw[very thick] (0,-1) -- node [right=-0.1] {$F_2$} (0,0) ;
\draw[very thick] (0,0) node [st1] {} -- node [right=-0.1] {$F_3$} (0,1);
\draw[thick] (0,1) node [st1] {} -- (1,0) node [st1] {} -- (0,-1) node [st1] {} -- (-1,0) node [st1] {} -- (0,1);\draw[very thick] (-1,0) -- node[below=-0.05] {$F_1$} (0,0);
\node at (0.5,0) {$K$};
\end{tikzpicture}
\qquad
\begin{tikzpicture}[scale=1.5,
st1/.style={circle,draw=black,fill=black,thick, minimum size=1mm, inner sep=0pt}
]
\draw[very thick] (0,-1) node [st1] {} -- node [right=-0.1] {$\bar{F}_2$} (0,-0.2) node [st1] {} ;
\draw[very thick] (0,0.2) node [st1] {} -- node [right=-0.1] {$\bar{F}_4$} (0,1) node [st1] {};
\draw[very thick] (0,-0.2) -- node [left=-0.1] {$\bar{F}_3$} (0,0.2);
\draw[thick] (0,1) node [st1] {} -- (1,0) node [st1] {} -- (0,-1) node [st1] {} -- (-1,0) node [st1] {} -- (0,1);
\draw[very thick] (-1,0) -- node [below=-0.05] {$\bar{F}_1$} (0,-0.2);
\draw[very thick] (1,0) -- node [below=-0.05] {$\bar{F}_5$} (0,0.2);
\end{tikzpicture}
\caption{Face regularity of meshes: the mesh on the left has interior faces $\calFhi=\left\{F_i\right\}_{i=1}^3$, each of which is regular to at least one element in the sense of Definition~\ref{def:face_regularity}, even though $F_2$ and $F_3$ fail to be regular with respect to the element $K$, since $F_2$ and $F_3$ are only proper subsets of the elemental face $F_2\cup F_3$. Since all boundary faces are also regular, the mesh on the left is face regular in the sense of Definition~\ref{def:face_regularity}, whereas the mesh on the right is not: the mesh face $\bar{F}_3$ fails to be regular with respect to any element of the mesh.}
\label{fig:face_regularity}
\end{center}	
\end{figure}
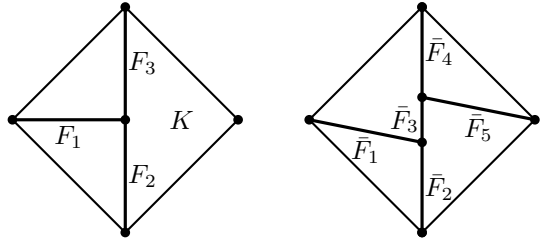

Let $\Om\subset \R^{\dim}$, $\dim\in\{2,3\}$, be a bounded polytopal domain with Lipschitz boundary $\DO$.
Let $\{\calT_h\}_{h>0}$ be a shape- and contact--regular sequence 
of simplicial meshes on $\Om$, as defined
in \cite[Definition 1.38]{DiPietro2012}.
For each element $\elm\in\calT_h$, let $h_\elm \coloneqq \diam \elm$, with $h = \max_{\elm\in\calT_h} h_\elm$ for each mesh $\calT_h$.
We define the faces of the mesh as in \cite[Definition~1.16]{DiPietro2012}, and we collect all interior and boundary faces in the sets $\calFhi$ and $\calFhb$, 
respectively, and let $\calFh \coloneqq \calF_h^i\cup\calF_h^b$ denote the skeleton of $\calT_h$. In particular, $F\in\calFhi$ if $F$ has positive $(\dim-1)$-dimensional Hausdorff measure and if $F=\partial K_1 \cap \partial K_2$ for two distinct mesh elements $K_1$ and $K_2$.
For an element $\elm \in \calT_h$, we denote $\calF(K)$ the set of faces of $K$, i.e.\ $E\in \calF(K)$ if $E$ is the closed convex hull of $d$ vertices of the simplex $K$. Note that on a mesh with hanging nodes, a mesh face may be a proper subset of an element face, see Fig.~\ref{fig:face_regularity}, hence the notions of mesh faces and element faces do not need to coincide. In this work, the meshes are allowed to have hanging nodes, provided that they satisfy the following notion of face regularity.

\begin{definition}\label{def:face_regularity}
A face $F\in\calFh$ is called \emph{regular} with respect to the element $K$ if $F\in \calF(K)$. We say that the mesh $\calT_h$ is \emph{face regular} if every face of $\calF_h$ is a regular face with respect to at least one element of $\calT_h$.
\end{definition}

Fig.~\ref{fig:face_regularity} illustrates the notion of face regularity with two examples. We remark that any matching mesh is face regular. On a face regular mesh, any boundary face is necessarily regular with respect to the element to which it belongs.
It appears that meshes of practical interest are most likely to be face regular, so this restriction is rather mild in practice.

 For integrable functions $\phi$ defined piecewise on either $\calT_h$ or $\calFh$, we use the convention
\[
\begin{aligned}
\int_\Om \phi \,\dd x = \sum_{\elm\in\calT_h} \int_\elm \phi \,\dd x,  & & &
\int_{\calF_h} \phi \,\dd s = \sum_{\face\in\calFh}\int_{\face}\phi \,\dd s.
\end{aligned}
\]

For the integer $k\geq 1$, we define the discontinuous finite element spaces 
$\Vh$ as the space of real-valued piecewise-polynomials of degree at most $k$ on $\calT_h$, 
and $\Sh$ the space of vector-valued piecewise-polynomials of degree at most $k+1$ on 
$\calT_h$. We define the mesh-dependent norm $\norm{\cdot}_{1,h}$ on $\Vh$ by
\begin{equation}\label{eq:dg_norm}
\begin{aligned}
\norm{v_h}_{1,h}^2 \coloneqq \sum_{\elm\in\calT_h} \norm{\nabla v_h}_{L^2(\elm)}^2
+ \sum_{\face\in\calFh} \frac{1}{h_F}\norm{\llb v_h \rrb}_{L^2(\face)}^2 & & & \forall\,v_h\in\Vh,
\end{aligned}
\end{equation}
where $h_F\coloneqq \diam \face$ for each face $\face\in\calFh$.

We shall also make use of the (local) Raviart--Thomas--N\'ed\'elec space \cite{Nedelec1980} defined by
\[
\LRTN\coloneqq \bcalP_{k}(\elm) \oplus \tilde{\calP}_k(\elm)\,\bx \subset \bcalP_{k+1}(\elm),
\]
where $\bcalP_k(\elm)$ is the space of vector-valued polynomials of 
degree at most $k$ on $\elm$, and $\tilde{P}_k(\elm)$ is the space of 
real-valued homogeneous polynomials of degree $k$ on $\elm$. 
We recall that $\btau_h\in  \LRTN$ is uniquely determined by the moments 
$\int_\elm \btau_h\cdot \bmu_h\,\dd x$ and $\int_E (\btau_h\cdot \bn_{E}) \,v_h\,\dd s$ 
for all $\bmu_h\in \bcalP_{k-1}(\elm)$ and $v_h\in \calP_k(E)$ for
each $E\in \calF(K)$, where $\bn_{E}$ denotes a unit normal vector of~$E$. We also recall that if all facial moments of $\btau_h$ vanish on an elemental face $E$, then $\btau_h\cdot \bn_{E}$ vanishes identically on $E$. 

For a face $\face\in\calFh$ belonging to an element 
$\elm_{\mathrm{ext}}$, we define the jump and average operators by
\begin{align*}
  \llb w \rrb|_F  &\coloneqq  \eval{w}{\elm_{\mathrm{ext}}} - \eval{w}{\elm_{\mathrm{int}}} ,
& \lla w \rra|_F &\coloneqq  \half\left( \eval{w}{\elm_{\mathrm{ext}}}  + \eval{w}{\elm_{\mathrm{int}}} \right) ,
&\text{if }\face\in\calF^i_h,
\\
\llb w \rrb|_F  &\coloneqq\eval{w}{\elm_{\mathrm{ext}}},
&  \lla w \rra|_F &\coloneqq  \eval{w}{\elm_{\mathrm{ext}}} ,
&\text{if }\face\in\calF_h^b,
\end{align*}
where $w$ is a sufficiently regular scalar or vector-valued function, and in the case where $\face\in\calFhi$, $\elm_{\mathrm{int}}$ is such that $\face= \p \elm_{\mathrm{ext}} \cap \p \elm_{\mathrm{int}}$.
Here, the labelling is chosen so that $\bn_F$ is outward pointing with respect to $\elm_{\mathrm{ext}}$ and inward pointing with respect to $\elm_{\mathrm{int}}$.
Let $\phi \in L^2(\calFh)$, then the lifting operators $\bo r_h \colon L^2(\calFh) \tends \Sh$ and $r_h\colon L^2(\calF_h)\tends \Vh$ are defined by
\begin{subequations}
\begin{align}
\int_\Om {\bo r_h}(\phi)\cdot \bsigma_h \, \dd x  & = \int_{\calFh} \phi  \lla \bsigma_h \cdot \bn_F \rra \,\dd s    & & \forall \,\bsigma_h \in \Sh,\label{eq:vector_lifting} \\
 \int_\Om r_h(\phi)\, v_h \,\dd x &= \int_{\calFhi} \phi \lla v_h \rra \,\dd s  & & \forall \,v_h\in\Vh. \label{eq:scalar_lifting}
\end{align}
\end{subequations}
For quantities $a$ and $b$, we write $a\lesssim b$ if and only if there is a positive constant $C$ such that $a\leq C b$, where $C$ is independent of the quantities of interest, such as the element sizes, but possibly dependent on the shape-regularity parameters and polynomial degrees.

\section{Stability of lifted gradients}\label{sec:stability}
We define the lifted gradient $G_h\colon \Vh \tends \Sh$ by
\begin{equation}\label{eq:lifted_gradient}
\begin{aligned}
G_h(v_h) = \nabla_h v_h - \bo r_h(\llb v_h \rrb)	 & & & \forall\, v_h\in \Vh,
\end{aligned}
\end{equation}
where $\nabla_h$ denotes the element-wise gradient operator. We note that $G_h$ is usually defined with a lifting using polynomial degrees $k$ or $k-1$, see for instance~\cite{DiPietro2012}. However, as we shall see, by increasing the polynomial degree of the lifting to $k+1$, we obtain the following key stability result.

\begin{theorem}\label{thm:stability}
Let $\{\calT_h\}_{h>0}$ denote a shape regular, contact regular and face regular sequence of simplicial meshes on $\Om$. 
Let the norm $\norm{\cdot}_{1,h}$ be defined by~\eqref{eq:dg_norm} 
and let the lifted gradient operator $G_h$ be defined by~\eqref{eq:lifted_gradient}. 
Then, we have 
\begin{equation}\label{eq:lifted_gradient_stability}
\begin{aligned}
\norm{u_h}_{1,h} \lesssim \norm{G_h(u_h)}_{L^2(\Om)} \lesssim \norm{u_h}_{1,h} & & &\forall\,u_h\in\Vh.
\end{aligned} 
\end{equation}
\end{theorem}
\begin{proof}
The upper bound $\norm{G_h(u_h)}_{L^2(\Om)}\lesssim \norm{u_h}_{1,h}$ is 
standard and we refer the reader to \cite[Sec.~4.3]{DiPietro2012} for a proof. 
To show the lower bound, consider an arbitrary $u_h\in \Vh$. Since $G_h(u_h)\in \Sh$, we have
\[
\norm{G_h(u_h)}_{L^2(\Om)} = \sup_{\btau_h\in \Sh\setminus\{0\}} \frac{\int_\Om G_h(u_h)\cdot \btau_h\,\dd x}{\norm{\btau_h}_{L^2(\Om)}},
\]
with the supremum being achieved by the choice $\btau_h=G_h(u_h)$. Therefore, to show \eqref{eq:lifted_gradient_stability}, it is sufficient to construct a $\btau_h\in\Sh$ such that
\begin{align}
\norm{u_h}_{1,h}^2 & \lesssim \int_{\Om} G_h(u_h)\cdot \btau_h\,\dd x, 	\label{eq:lower_bound} \\
\norm{\btau_h}_{L^2(\Om)} & \lesssim \norm{u_h}_{1,h}. \label{eq:upper_bound}
\end{align}
Let $\btau_\elm\in \LRTN$ be defined by
\begin{subequations}\label{eq:tau_h_definition}
\begin{align}
\int_\elm \btau_\elm \cdot \bmu_h \,\dd x & = \int_{\elm} \nabla u_h \,\cdot \bmu_h \,\dd x\quad   \forall\, \bmu_h\in \bcalP_{k-1}(\elm), \label{eq:tau_h_definition1}
\\
\int_E \left(\btau_\elm\cdot \bn_{E}\right) v_h\,\dd s & =
\begin{cases}
- \int_E \frac{1}{h_E} \llb u_h \rrb \, v_h \,\dd s  & \text{if } E\in \calFh, \\
 0 & \text{if } E\notin \calFh,
\end{cases}
\label{eq:tau_h_definition2}
\end{align}
\end{subequations}
where \eqref{eq:tau_h_definition2} holds for all $v_h\in \calP_k(E)$, for each element face $E\in\calF(K)$. In particular, if the element face $E\in\calFh$, i.e. $E$ is also a mesh face, then we require that $\bn_E$ agrees with the choice of unit normal used to define the jump and average operators. If $E\notin \calFh$, then $\btau_\elm\cdot \bn_{E}$ vanishes identically on $E$, and the orientation of $\bn_E$ on the left-hand side of \eqref{eq:tau_h_definition2} does not matter.
The global vector field $\btau_h\in \Sh $ is defined element-wise by $\btau_h|_{\elm} = \btau_\elm$.

Since the mesh $\calT_h$ is assumed to be face regular, for every $F\in\calFh$ there exists an element $K\in\calT_h$ and an elemental face $E\in \calF(K)$ such that $E=F$; then $E$ satisfies the first condition in \eqref{eq:tau_h_definition2}.
Therefore, the facts that $\lla \btau_h\cdot \bn_F\rra|_F$ and $\llb u_h \rrb|_F$ both belong to $\calP_k(F)$ together with \eqref{eq:tau_h_definition2} imply that for each $F\in\calFh$, one of only three situations may arise:
\begin{enumerate}
\item $F$ is a boundary face and hence $F\in \calF(K)$. In this case, we have $\lla \btau_h\cdot \bn_F\rra |_F = -h_F^{-1} \llb u_h\rrb|_F$.
\item $F$ is an interior face which is regular with respect to both elements to which it belongs. In this case, we have $\lla \btau_h\cdot \bn_F\rra |_F = -h_F^{-1} \llb u_h\rrb|_F$.
\item $F$ is an interior face which is regular with respect to only one of the elements to which it belongs. In this case, we have $\lla \btau_h\cdot \bn_F\rra|_F = -\half h_F^{-1} \llb u_h\rrb|_F$, since $\eval{\btau_{h}}{K^\prime}\cdot \bn_F\equiv 0$ for the element $K^{\prime}$ with respect to which $F$ is not regular.
\end{enumerate}

Therefore, since $\btau_h\in \Sh$, the definition of the lifting operator in \eqref{eq:vector_lifting} implies that
\begin{equation*}
\begin{split}
\int_\Om G_h(u_h) \cdot \btau_h\,\dd x 
& = \sum_{\elm\in\calT_h} \int_\elm \nabla u_h \cdot \btau_h \,\dd x - \sum_{\face\in\calFh}\int_F \{\btau_h\cdot \bn_F\}\,\llb u_h \rrb\,\dd s \\
& \ge  \sum_{\elm\in\calT_h} \norm{\nabla u_h}_{L^2(\elm)}^2 + \frac{1}{2} \sum_{\face\in\calFh}\frac{1}{h_F} \norm{\llb u_h \rrb}_{L^2(\face)}^2 \\
& \ge \frac{1}{2} \norm{u_h}_{1,h}^2,
\end{split}
\end{equation*}
where the second line follows from \eqref{eq:tau_h_definition} and 
from the fact that $\eval{\nabla u_h}{\elm}\in \bcalP_{k-1}(\elm)$ for each $\elm\in\calT_h$.
Hence  \eqref{eq:lower_bound} is satisfied,
and we now verify \eqref{eq:upper_bound}.
A classical scaling argument using the Piola transformation \cite[p.~59]{Boffi2013} yields
\begin{multline*}
\norm{\btau_h}_{L^2(\elm)} \lesssim \sup_{\bmu_h\in \bcalP_{k-1}(\elm)\setminus\{0\}} \frac{\int_\elm \btau_h\cdot \bmu_h\,\dd x }{\norm{\bmu_h}_{L^2(\elm)}}\\
 + \sum_{E\in\calF(K)}\sup_{v_h\in \calP_k(E)\setminus\{0\}} \frac{h^{1/2}_E\int_E (\btau_h\cdot \bn_E) v_h \,\dd s }{\norm{v_h}_{L^2(E)}} \quad\forall\,\elm\in\calT_h.
\end{multline*}
Therefore, it follows from \eqref{eq:tau_h_definition} that, for each $K\in\calT_h$,
\begin{equation}\label{eq:upper_bound_element}
\norm{\btau_h}_{L^2(\elm)}^2 \lesssim \norm{\nabla u_h}_{L^2(\elm)}^2 +\sum_{F\in \calF(K)\cap \calFh} h_F \norm{ h_F^{-1} \llb u_h \rrb}_{L^2(F)}^2.
\end{equation}
Summing \eqref{eq:upper_bound_element} over all elements therefore implies \eqref{eq:upper_bound}. \qed
\end{proof}

\section{Counterexample to stability for equal-order liftings}\label{sec:counter}
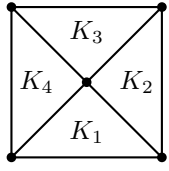
\begin{figure}[t]
\begin{center}
\begin{tikzpicture}[scale=1,
st1/.style={circle,draw=black,fill=black,thick, minimum size=1mm, inner sep=0pt}
]
\draw[thick] (-1,-1) node [st1] {} --(1,-1) node [st1] {} --(1,1) node [st1] {} --(-1,1) node [st1] {} --(-1,-1) ;
\draw[thick] (-1,1) -- (1,-1);
\draw[thick] (-1,-1) -- (0,0) node [st1] {} -- (1,1);
\node at (0,-0.6667) {$K_1$};
\node at (0.6667,0) {$K_2$};
\node at (0,0.6667) {$K_3$};
\node at (-0.6667,0) {$K_4$};
\end{tikzpicture}
\caption{Counterexample of Sec.~\ref{sec:counter}: the domain $\Om=(-1,1)^2$ and the criss-cross mesh $\calT_h$ considered in the example.}
\label{CounterExamplePic}
\end{center}
\end{figure}

Theorem~\ref{thm:stability} shows the stability of the lifted gradient operator $G_h$ provided that the lifting operator $\bo r_h$ has polynomial degree $k+1$. In this section, we verify by means of a counterexample that the stability estimate does not generally hold for lower-order liftings, including in particular the case of equal-order liftings, which are commonly used in practice; our example  simplifies a similar counterexample in \cite{Brezzi1997}.

{\em Example.} Let $\Omega = (-1,1)^2$, and consider the finite element space $V_{h,k}$ defined on a criss-cross mesh with four triangles, as depicted in Fig.~\ref{CounterExamplePic}, using piecewise linear polynomials, i.e.\ $k=1$. Let $u_h\in V_{h,1}$ be the piecewise linear function defined by
\begin{align*}
u_h |_{K_1} &= y+\frac{2}{3},  & 
u_h |_{K_2} &= x-\frac{2}{3}, \\
u_h |_{K_3} &= -y+\frac{2}{3}, &
u_h |_{K_4} &= -x - \frac{2}{3}.
\end{align*}
Direct calculations show that $\{u_h\}|_F\equiv 0$ on all interior faces $F\in\calFhi$, and that $\int_K u_h \,\dd x =0$ for all elements $K\in\calT_h$. Consequently, if the lifting operator $\bo{\tilde{r}}_h$ is defined in \eqref{eq:vector_lifting} with the polynomial degree $k+1$ replaced by $k$, and if $\tilde{G}_h(u_h)\coloneqq \nabla_h u_h - \bo{\tilde{r}}_h( \llb u_h \rrb )$ denotes the equal-order lifted gradient, then we have for all $\btau_h \in \bm{\Sigma}_{h,1}$,
\begin{align*}
\int_\Om \tilde{G}_h(u_h) \cdot \btau_h\,\dd x 
& = \sum_{\elm\in\calT_h} \int_\elm \nabla_h u_h \cdot \btau_h \,\dd x - \sum_{\face\in\calFh}\int_F \{\btau_h\cdot \bn_F\}\,\llb u_h \rrb\,\dd s \\
&= -\sum_{\elm\in\calT_h} \int_\elm u_h (\nabla_h \cdot \btau_h)\, \dd x +\sum_{\face\in\calFh^i}\int_F \{u_h\} \llb \btau_h\cdot \bn_F\rrb\, \dd s =0.
\end{align*}
Since $\tilde{G}_h(u_h)\in\bm{\Sigma}_{h,1}$, we deduce that $\tilde{G}_h(u_h)=0$, and thus it is found that no bound of the form $\norm{u_h}_{1,h}\lesssim \norm{\tilde{G}_h(u_h)}_{L^2(\Om)}$ is possible.\qed

\section{A modified LDG method without penalty parameters}\label{sec:LDG}
As an application of Theorem~\ref{thm:stability}, consider the  discretization of the homogeneous Dirichlet boundary-value problem of the Poisson equation by a modified LDG method \cite{Castillo2000,Cockburn1998} as follows. For $f\in L^2(\Om)$, let $u\in H^1_0(\Om)$ be the unique solution of 
\begin{equation}\label{eq:poisson}
\int_\Om \nabla u\cdot \nabla v\,\dd x = \int_\Om f\,v\,\dd x \quad \forall\, v \in H_0^1(\Om).	
\end{equation}
Let the bilinear form $a_h \colon \Vh \times \Vh \tends \R$ be defined by
\begin{equation}
\begin{aligned}
a_h(u_h,v_h) = \int_\Om G_h(u_h)\cdot G_h(v_h) \,\dd x & & & \forall\, u_h,\,v_h\in\Vh,
\end{aligned}
\end{equation}
where the lifted gradient operator $G_h$ was defined in \eqref{eq:lifted_gradient}. The bilinear form $a_h(\cdot, \cdot)$ defines a modified LDG method for \eqref{eq:poisson}: find $u_h\in\Vh$ such that
\begin{equation}\label{eq:mldg}
\begin{aligned}
	a_h(u_h,v_h) = \int_\Om f\,v_h\,\dd x & & & \forall v_h\in\Vh.
\end{aligned}
\end{equation}
It follows from Theorem~\ref{thm:stability} that $a_h(\cdot, \cdot)$ is uniformly stable with respect to the norm $\norm{\cdot}_{1,h}$, and thus \eqref{eq:mldg} is well-posed for each $h$. Moreoever, the discrete Poincar\'e inequality \cite{DiPietro2012} implies that $\norm{u_h}_{1,h}\lesssim \norm{f}_{L^2(\Om)}$ for all $h$, so that the numerical solutions $u_h$ are uniformly bounded with respect to the mesh-dependent norms $\norm{\cdot}_{1,h}$. The \emph{a priori} error analysis for the numerical method defined by \eqref{eq:mldg} may be developed following the frameworks of \cite{Castillo2000,DiPietro2012,Gudi2010}, although for reasons of space we do not present the arguments here.

An interesting feature of the modified LDG method~\eqref{eq:mldg} is that it does not require any additional stabilization, such as added penalty terms of the form $\int_{\calFh} \tfrac{\sigma_F}{h_F} {\llb u_h \rrb}{\llb v_h \rrb}\, \dd s$ for some user-defined parameter $\sigma_F$.
The absence of such penalty terms enables us to show the following discrete conservation property. We define the lifted divergence $D_h\colon \Sh\tends \Vh$ by
\begin{equation}
\begin{aligned}
D_h(\bsigma_h) = \Div_h \bsigma_h - r_h(\llb \bsigma_h\cdot \bn_F \rrb), & & & \bsigma_h\in \Sh,
\end{aligned}
\end{equation}
where $\Div_h$ denotes the element-wise divergence operator, and where $r_h$ is the scalar lifting operator defined in \eqref{eq:scalar_lifting}. We note that we have the integration-by-parts identity
\begin{equation}
\begin{aligned}
\int_\Om \bsigma_h \cdot G_h(v_h)\,\dd x = - \int_\Om D_h(\bsigma_h) \,v_h \,\dd x  & & & \forall\, v_h\in\Vh,\,\bsigma_h\in\Sh,
\end{aligned}
\end{equation}
which should be compared with the analogous continuous identity between the spaces $H^1_0(\Om)$ and $H(\Div,\Om)$.
Therefore, the numerical scheme \eqref{eq:mldg} can be equivalently expressed in the strong form
\begin{equation}
-\int_\Om D_h(G_h(u_h))\, v_h\,\dd x = \int_\Om f\, v_h\,\dd x,
\end{equation}
which implies that the numerical solution $u_h\in\Vh$ solves
\begin{equation}
-D_h(G_h(u_h))=\Pi_{h}^k f,
\end{equation}
in the pointwise sense on each element $\elm$, where $\Pi_h^k f$ denotes the element-wise $L^2$-projection of $f$ into $\Vh$.
Although we have shown here how the lifted gradient operator $G_h$ of degree $k+1$ may be used to achieve a stable discretization of the Poisson equation, it is by no means restricted to this model problem, as the lifted gradients may be used to discretize the second-order terms of more general differential operators.

\section{Conclusions}\label{sec:conclusions}
In this article, we studied an intrinsically stable modified LDG method without additional parameter dependent penalization. For this, we showed that increasing the degree of the lifting operator by one order leads to stability of the discrete gradient operator on face regular meshes with hanging nodes.

\subsection*{Acknowledgement}
The work of the third author
was partially supported by the NSF grant DMS--1417980
and the Alfred Sloan Foundation.


\end{document}